\documentclass[12pt]{article}
\usepackage{my} 
\usepackage{float}
\usepackage{enumerate}

\title{The maximum agreement subtree problem}

\author{
  Daniel M. Martin \thanks{Supported by CNPq (Processo: 475064/2010-0)} \\
  \small  Centro de Matem\'{a}tica, Computa\c{c}\~{a}o e Cogni\c{c}\~{a}o \\
  \small Universidade Federal do ABC, \\
  \small Santo Andr\'{e}, SP 09210-170 Brasil. \\
  \small \texttt{daniel.martin@ufabc.edu.br} \and Bhalchandra D. Thatte
  \thanks{Supported by CNPq (Processo: 151782/2010-5) and partially
    supported by Project MaCLinC (Mathematics, computation, language and
    the brain) at USP}\\
  \small Instituto de Matem\'atica e Estat\'istica,\\
  \small Universidade de S\~{a}o Paulo, \\
  \small SP 05508-090 Brasil.\\
  \small \texttt{bdthatte@gmail.com} }
\date{\today \\[2mm]
  \small Mathematics Subject Classifications: Primary: Secondary: }

\usepackage{amsmath,amsthm}
\usepackage{latexsym,amssymb}
\usepackage{graphicx, epsfig}
\usepackage{mathtools}
\usepackage{caption, subcaption}

\theoremstyle{plain}
\newtheorem{thm}{Theorem}

\newtheorem{cor}[thm]{Corollary}
\newtheorem{lem}[thm]{Lemma}
\newtheorem{prop}[thm]{Proposition}
\newtheorem{conj}[thm]{Conjecture}
\newtheorem{claim}[thm]{Claim}

\theoremstyle{definition}

\newtheorem{problem}[thm]{Problem}

\theoremstyle{remark}
\newtheorem{rem}[thm]{Remark}

\newcommand{\mcb}{\mathcal{B}}
\newcommand{\mcc}{\mathcal{C}}

\newcommand{\mct}{\mathcal{T}}

\usepackage[numbers]{natbib} 

\newcommand{\ta}[1]{T_1^{(#1)}}
\newcommand{\tb}[1]{T_2^{(#1)}}
\newcommand{\mast}{{\mathrm{mast}}}

\begin{document}
\maketitle
\begin{abstract}
  In this paper we investigate an extremal problem on binary
  phylogenetic trees. Given two such trees $T_1$ and~$T_2$, both with
  leaf-set $\{1,2,\dots,n\}$, we are interested in the size of the
  largest subset $S \subseteq \{1,2,\dots,n\}$ of leaves in a common
  subtree of~$T_1$ and~$T_2$. We show that any two binary phylogenetic
  trees have a common subtree on $\Omega(\sqrt{\log{n}})$ leaves, thus
  improving on the previously known bound of $\Omega(\log\log n)$ due to
  M.~Steel and L.~Szekely. To achieve this improved bound, we first
  consider two special cases of the problem: when one of the trees is
  balanced or a caterpillar, we show that the largest common subtree has
  $\Omega(\log n)$ leaves. We then handle the general case by proving
  and applying a Ramsey-type result: that every binary tree contains
  either a large balanced subtree or a large caterpillar.  We also show
  that there are constants $c, \alpha > 0$ such that, when both trees
  are balanced, they have a common subtree on $c n^\alpha$ leaves. We
  conjecture that it is possible to take $\alpha = 1/2$ in the unrooted
  case, and both $c = 1$ and $\alpha = 1/2$ in the rooted case.
\end{abstract}

\section{Preliminaries}\label{sec-pre}

All trees considered in this paper are binary. Although we mainly talk
about rooted trees, we introduce the problem in terms of unrooted trees
to be consistent with earlier papers on the subject.

\subsection {Unrooted phylogenetic trees}

A \emph{phylogenetic tree} is a binary, unrooted tree in which the
leaves are labelled bijectively with elements from a finite set. All
internal vertices of a phylogenetic tree have degree~$3$. For such a
tree~$T$, the set of vertices is denoted by~$V(T)$, the set of edges
by~$E(T)$, and the set of leaves by~$L(T)$.

In phylogenetics, it is common to consider isomorphism between trees in
a more restricted sense. We say that trees $T_1$ and $T_2$ are
\emph{isomorphic} (and write $T_1 \cong T_2$) if there is a bijection
$\varphi: V(T_1) \rightarrow V(T_2)$ such that
\begin{enumerate}[(i)]
\item $\{\varphi(u),\varphi(v)\} \in E(T_2) \iff \{u,v\} \in E(T_1)$,
\item $\varphi(i) = i$ for all leaves $i \in L(T_1)$.
\end{enumerate} Observe that, while this notion of isomorphism also
works for non-binary trees, there can be no isomorphism between trees
that have distinct leaf-sets.

\begin{figure}[H]
  \centering
  \includegraphics{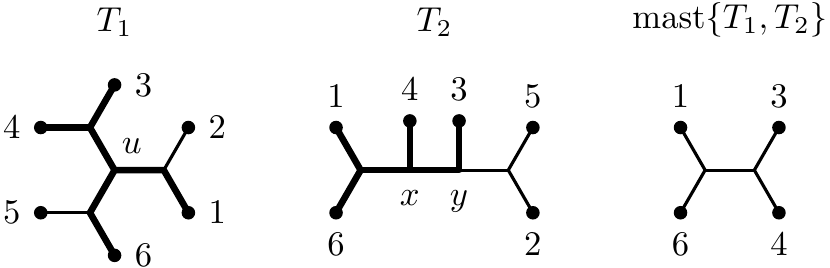}
  \caption{}
  \label{fig:example_mast}
\end{figure}

For a subset $X \subseteq L(T)$, define the {\em restriction} of $T$ to
$X$ to be the phylogenetic tree $T|X$ with leaf-set $L(T|X) = X$, and
the property that there exists an isomorphism (in the sense of the
previous paragraph) from a subdivision of $T|X$ to the unique minimal
connected subgraph of $T$ containing~$X$. We loosely call the tree $T|X$
a \emph{subtree} of $T$. Given trees $T_1$ and $T_2$, if $X$ is a subset
of~$L(T_1) \cap L(T_2)$ of maximum cardinality with the property that
$T_1|X \cong T_2|X$, we say that $T_1|X$ (or $T_2|X$) is a {\em maximum
agreement subtree} of $T_1$ and $T_2$. We also define the
parameter~$\mast\{T_1,T_2\}$ as
\[ \mast\{T_1,T_2\} \coloneqq \max \big\{ |X| \colon X \subseteq L(T_1)
\cap L(T_2), T_1|X \cong T_2|X \big\}.
\] 
Figure~\ref{fig:example_mast} illustrates trees $T_1$ and~$T_2$ and
one of their maximum agreement subtrees.

Now let \[\mast(n) \coloneqq \min \big\{ \mast\{T_1, T_2\} \colon
L(T_1) = L(T_2) = \{1,2,\dots,n\} \big\}.\] It was shown by
\citeauthor*{kkm1992} in \cite{kkm1992} that\footnote{Logarithms in this
paper are always taken with base~$2$.}
\begin{equation}
  \label{eq:kkm} c_1 (\log \log n)^{1/2} \leq \mast(n) \leq c_2 \log n,
\end{equation}
for some positive constants $c_1$ and $c_2$. To see that the upper bound
in~\eqref{eq:kkm} is tight, consider the case when $T_1$ is a
caterpillar with $n$ leaves and $T_2$ is a balanced tree of height $\log
n$ (to be defined in Section~\ref{sec:onebalanced}): any common subtree
must be a caterpillar, and there is no caterpillar of length more than
$2\log n$ in~$T_2$. The lower bound in~\eqref{eq:kkm} was improved by
\citet{ss2009} who showed that $\mast(n) \geq c \log \log n$ for a
positive constant $c$. In fact, in a remark following Theorem~$1$ in
their paper, they mention that a more explicit bound of $\frac{1}{4}\log
\log (n-1)$ may be derived, and suggest that a much stronger lower bound
of $c\log n$ might hold, for some positive constant $c$.

\begin{problem}\label{c} Is there a constant $c > 0$ such that any two
  phylogenetic trees $T_1$ and $T_2$ with $L(T_1) = L(T_2) =
  \{1,2,\dots,n\}$ have a maximum agreement subtree on at least $c \log
  n$ leaves?
\end{problem}

One of the goals of this paper is to further improve the lower bound
in~\eqref{eq:kkm}. In particular, in Theorem~\ref{thm:root_log}, we show
that any two phylogenetic trees $T_1$ and $T_2$ with leaf-set
$\{1,2,\dots,n\}$ have an agreement subtree on $\Omega(\sqrt{\log n})$
leaves.

\subsection {Rooted phylogenetic trees}

Next, we develop some terminology for rooted phylogenetic trees. In a
rooted phylogenetic tree~$T$, for $|L(T)| > 1$, all internal nodes have
degree~$3$ except the root, which has degree~$2$. If $|L(T)| = 1$,
then~$T$ has exactly one vertex which is both its only leaf and
root. Let us denote the root of a tree~$T$ by~$\rho(T)$. For a vertex~$u
\in V(T)$, we denote by~$T^{u}$ the subtree of~$T$ rooted at~$u$
containing all descendants of~$u$ in $T$. We denote by $\ell(u)$ and
$r(u)$, respectively, the left and the right children of an internal
vertex $u$ of a rooted tree.

For a rooted tree $T$ and vertices $x$ and $y$ in $V(T)$, we define $x
\wedge y$ to be the \emph{most recent common ancestor} of $x$ and
$y$. Since $\wedge$ is associative and commutative, we may define the
most recent common ancestor of a set $X = \{x_1, x_2, \dots, x_k\}
\subseteq V(T)$ to be
$$
\bigwedge X \coloneqq x_1 \wedge x_2 \wedge \cdots \wedge x_k.
$$ 
Given rooted trees $S$ and $T$, we say that $S$ is a subtree of $T$ (and
write $S \preceq T$) if there exists an injective map $f:
V(S)\rightarrow V(T)$ satisfying
\begin{enumerate}
\item $f(x) = x$ for all $x \in L(S)$
\item $f(x \wedge y) = f(x) \wedge f(y)$ for all $x,y \in V(S)$
\end{enumerate}
We say that rooted phylogenetic trees $T_1$ and $T_2$ are
\emph{isomorphic} if $T_1 \preceq T_2$ and $T_2 \preceq T_1$.

Now we define the \emph{restriction} of $T$ to the set of leaves $X$ as
the unique binary, rooted phylogenetic tree $T|X$ having leaf-set $X$
and satisfying $T|X \preceq T$. The rooted maximum agreement subtree and
$\mast\{\cdot\}$ are defined as in the unrooted case.

Proposition~\ref{lemma:join} allows us to recursively construct
agreement subtrees (but not necessarily maximum agreement subtrees) of
rooted trees.

\begin{prop}\label{lemma:join} Let $T_1$ and $T_2$ be rooted, binary
trees, with roots $u \coloneqq \rho(T_1)$ and $v \coloneqq \rho(T_2)$,
respectively. If $S_\ell \preceq \ta{\ell(u)}, \tb{\ell(v)}$ and $S_r
\preceq \ta{r(u)}, \tb{r(v)}$, then the tree $S$ with left and right
subtrees $S_\ell$ and $S_r$, respectively, is such that $S \preceq T_1,
T_2$.
\end{prop}

\begin{proof} Let $f_\ell$ and $f_r$ be maps that realize $S_\ell
  \preceq \ta{\ell(u)}$ and $S_r \preceq \ta{r(u)}$, respectively. Then
  $S \preceq T_1$ is realized by defining a map $f: V(S) \rightarrow
  V(T_1)$ so that
  \[ f(x) \coloneqq
  \begin{cases} f_\ell(x) & \text{ if } x \in V(S_\ell) \\
    f_r(x) & \text{ if } x \in V(S_r) \\ \rho(T_1) & \text{ if } x =
    \rho(S).
  \end{cases}
  \] A similarly constructed map shows that $S \preceq T_2$.
\end{proof}

We denote the tree $S$ constructed from $S_\ell$ and $S_r$ as in the
above lemma by $S_\ell \circ S_r$.

To prove the results in this paper, we first obtain agreement subtrees
of rooted trees (constructed by rooting the given unrooted trees
suitably), and then agreement subtrees of unrooted trees by ignoring the
roots.

\section{When one of the trees is balanced}
\label{sec:onebalanced}

A rooted, phylogenetic tree is \emph{balanced} if all leaves are at
the same distance from the root. For unrooted trees, the definition is
analogous. We first define the center of a graph~$G$ as the set of
vertices $u \in V(G)$ for which
\[
\max \{ d_G(u,v) \colon v \in V(G) \}
\]
is minimum, where $d_G(u,v)$ is the length of a shortest path from~$u$
to~$v$ in~$G$. The center of a tree contains either a single vertex or
two adjacent vertices. For example, in Figure~\ref{fig:example_mast},
$T_1$ has center $\{ u \}$, while $T_2$ has center $\{x, y\}$. We say
that a phylogenetic tree is \emph{balanced} if all leaves are at the
same distance from the center of the tree. For example, the tree $T_1$ in
Figure~\ref{fig:example_mast} is balanced. Note that a rooted
balanced tree has $2^m$ leaves for some $m \geq 0$, while an unrooted
balanced tree has either $2^m$ or $3 \cdot 2^{m}$ leaves for some $m
\geq 0$. We refer to Section~\ref{sec:unrooted_case} for more
explanation.

In this section, we solve Problem~\ref{c} when one of the trees is
balanced and binary.

\subsection {The rooted case}

We first consider the case of rooted trees, thus assuming that one tree
is a rooted, balanced tree of height $m$ and another is a general
rooted, binary tree. In Theorem~\ref{t1} we prove a result for the
unrooted case by appropriately rooting the trees and applying
Lemma~\ref{lemma:match1}.

In what follows, for nodes $x \in V(T_1)$ and~$y \in V(T_2)$, let
$t_{xy}$ be the number of elements in the set $L(T_1^x) \cap L(T_2^y)$.

\begin{lem}
  \label{lemma:match1} Suppose $T_1$ is a rooted, balanced, binary tree
on a leaf-set of cardinality $2^m$, and $T_2$ is an arbitrary rooted,
binary tree on $t > 0$ leaves with $L(T_2) \subseteq L(T_1)$. Then for
all $\delta \in (0,1/2)$, the two trees have a maximum agreement subtree
that has at least
  $$
  \displaystyle \frac{m\log(1-\delta) + \log{t}}{1-\log{\delta}}
  $$
  leaves.
\end{lem}

\begin{proof} Let $g(m,t)$ be the minimum value of $\mast\{T_1,T_2\}$
(over all choices of $T_1$ and $T_2$), where $T_1$ and $T_2$ are as in
the statement of the lemma. Observe that $g(m,t)$ is a monotonically
non-decreasing function of $t$. We show by induction on $m$ that

  \begin{equation}\label{eq:match1} g(m,t) \geq \frac{m\log(1-\delta) +
\log{t}}{1-\log{\delta}}.
  \end{equation}

  \

  \noindent {\em Base case:} If $m = 0$, then $g(m,t) = 1$ and the
right-hand-side of~\eqref{eq:match1} is 0. So we may assume that $m >
0$.

  \

  \noindent {\em Induction step:} Let $u \coloneqq \rho(T_1)$ and $v
\coloneqq \rho(T_2)$. Observe that
  \begin{equation}
    \label{eq:common_leaves} t = t_{uv} = t_{\ell(u)\ell(v)} +
    t_{r(u)r(v)} + t_{\ell(u)r(v)} + t_{r(u)\ell(v)}.
  \end{equation} Without loss of generality, assume that
  $t_{\ell(u)\ell(v)} + t_{r(u)r(v)} \geq t_{\ell(u)r(v)} +
  t_{r(u)\ell(v)}$, and $t_{r(u)r(v)} \geq t_{\ell(u)\ell(v)}$. Therefore,
  by~\eqref{eq:common_leaves}, we have
  \begin{equation}
    \label{eq:case1a} t_{r(u)r(v)} \geq \lceil t/4 \rceil.
  \end{equation}

  \

  \noindent {\em Case 1:} $t_{\ell(u)\ell(v)} > 0$.
  
  In this case, we take a maximum agreement subtree $S_\ell$ of
  $\ta{\ell(u)}$ and $\tb{\ell(v)}$, and a maximum agreement subtree
  $S_r$ of $\ta{r(u)}$ and $\tb{r(v)}$. We then construct $S = S_\ell
  \circ S_r$, which by Proposition~\ref{lemma:join} is an agreement
  subtree of $T_1$ and $T_2$. Now
  \[ g(m,t) \geq |L(S)| = |L(S_\ell)| + |L(S_r)| \geq 1 + g(m-1,
  t_{r(u)r(v)}).
  \] By~\eqref{eq:case1a}, we have
  $$
  g(m,t) \geq 1 + g(m-1,\lceil t/4 \rceil).
  $$
  Therefore, applying the induction hypothesis, we obtain
  \begin{eqnarray*}
    \label{eq:case1b} g(m,t) & \geq & 1 +
    \frac{(m-1)\log{(1-\delta)}+\log{(t/4)}} {1-\log{\delta}} \nonumber \\ &
    = & \frac{m\log{(1-\delta)}+\log{t}}{1-\log{\delta}} +
    \frac{-1-\log{(1-\delta)}-\log \delta}{1-\log{\delta}} \nonumber \\ & >
    & \frac{m\log{(1-\delta)}+\log{t}}{1-\log{\delta}}.
  \end{eqnarray*}
  
  \

  \noindent {\em Case 2:} $t_{\ell(u)\ell(v)} = 0$.

  From now on, we assume that $t_{r(u)\ell(v)}$ is non-zero; otherwise,
  together with the assumption of this case, it would contradict the
  hypothesis that $L(T_2) \subseteq L(T_1)$.

  \noindent {\em Subcase 2.1:} $t_{\ell(u)r(v)} + t_{r(u)\ell(v)} \geq
  \delta t$ and $t_{\ell(u)r(v)} > 0$. Also, without loss of generality,
  assume that $t_{r(u)\ell(v)} \geq t_{\ell(u)r(v)}$, so that
  $t_{r(u)\ell(v)} \geq \delta t/2$.

  In this case, we construct $S \coloneqq S_\ell \circ S_r$, where
  $S_\ell$ is a rooted maximum agreement subtree of $\ta{\ell(u)}$ and
  $\tb{r(v)}$, and $S_r$ is a rooted maximum agreement subtree of
  $\ta{r(u)}$ and $\tb{\ell(v)}$. Therefore,
  \begin{eqnarray*}
    \label{eq:case2.1} g(m,t) & \geq & |L(S_\ell)| + |L(S_r)| \nonumber
    \\ & \geq & 1 + g(m-1,\lceil \delta t/2 \rceil) \nonumber \\ & \geq & 1
    + \frac{(m-1)\log{(1-\delta)}+\log{(\delta t/2)}} {1-\log{\delta}}
    \nonumber \\ & \geq & \frac{m\log{(1-\delta)}+\log{t}}{1-\log{\delta}} +
    \frac{-\log{(1-\delta)}}{1-\log{\delta}} \nonumber \\ & > &
    \frac{m\log{(1-\delta)}+\log{t}}{1-\log{\delta}}.
  \end{eqnarray*}

  \

  \noindent {\em Subcase 2.2:} $t_{\ell(u)r(v)} + t_{r(u)\ell(v)} <
  \delta t$. In this case, $t_{r(u)r(v)} = t - t_{\ell(u)\ell(v)} -
  t_{\ell(u)r(v)} - t_{r(u)\ell(v)} > (1-\delta)t$.

  Let $S$ be a rooted maximum agreement subtree of $\ta{r(u)}$ and
  $\tb{r(v)}$. We have
  \begin{eqnarray*}
    \label{eq:case2.2} g(m,t) & \geq & |L(S)| \nonumber \\ & = &
    g(m-1,\lceil (1-\delta)t \rceil) \nonumber \\ & \geq &
    \frac{(m-1)\log{(1-\delta)} + \log{(1-\delta)t}}{1-\log{\delta}}
    \nonumber \\ & = & \frac{m\log{(1-\delta)}+\log{t}}{1-\log{\delta}}
  \end{eqnarray*}

  \

  \noindent {\em Subcase 2.3:} $t_{\ell(u)r(v)} = 0$.

  Let $S$ be a rooted maximum agreement subtree of $\ta{r(u)}$ and
  $\tb{v}$. Therefore,
  \begin{eqnarray*}
    \label{eq:case2.3} g(m,t) & \geq & |L(S)| \nonumber \\ & = &
    g(m-1,t) \nonumber \\ & \geq &
    \frac{(m-1)\log{(1-\delta)}+\log{t}}{1-\log{\delta}} \nonumber \\ & = &
    \frac{m\log{(1-\delta)}+\log{t}}{1-\log{\delta}} +
    \frac{-\log{(1-\delta)}}{1-\log{\delta}} \nonumber \\ & > &
    \frac{m\log{(1-\delta)}+\log{t}}{1-\log{\delta}}.
  \end{eqnarray*}
  With this we complete all subcases of the induction
  step. Therefore~\eqref{eq:match1} holds and the lemma is proved.
\end{proof}

An immediate consequence of the above lemma is the following corollary.

\begin{cor}
  \label{cor:match1} For $\delta \in (0,1/2)$, set
  \begin{equation} \displaystyle \alpha \coloneqq \frac{1 +
      \log{(1-\delta)}}{1 - \log{\delta}}.\label{eq:alpha}
  \end{equation} If $T_1$ is a rooted, balanced, binary tree on a
  leaf-set of cardinality $2^m$, and $T_2$ is an arbitrary rooted, binary
  tree on $2^m$ leaves such that $L(T_2) = L(T_1)$, then $T_1$ and $T_2$
  have a maximum agreement subtree on at least $\alpha m$ leaves.
\end{cor}

\begin{rem}
  \label{rem:alpha} We found numerically that the maximum value of
  $\alpha$ is approximately $0.2055$ obtained when $\delta$ is
  approximately $0.1705$.
\end{rem}

\begin{rem} An algorithm to construct an agreement tree (though not
  necessarily a maximum agreement subtree) is implicit in the proof of
  Lemma~\ref{lemma:match1}. Observe that in each of Case 1 and Subcase
  2.1, we may take the agreement subtree $S_\ell$ to be a tree with a
  single leaf (e.g. any leaf from $L(\ta{\ell(u)})\cap L(\tb{\ell(v)})$
  in Case 1 and any leaf from $L(\ta{\ell(u)})\cap L(\tb{r(v)})$ in
  Subcase 2.1). Such a choice gives us an agreement subtree that is a
  caterpillar of length $\alpha m$ in Corollary~\ref{cor:match1}.
\end{rem}

We explicitly describe a recursive algorithm, which we call
\texttt{Match1}. The algorithm takes as input two rooted, binary trees,
the first one being balanced, and returns a set of leaves in a common
subtree that is a caterpillar. As in Lemma~\ref{lemma:match1}, algorithm
\texttt{Match1} also depends on a parameter $\delta \in (0,1/2)$.
\medskip

\noindent \texttt{Algorithm Match1($\ta{u}$, $\tb{v}$)}
\begin{enumerate}[1:]

\item \label{m1a}\texttt{if $|L(\ta{u})| = 1$ or $|L(\ta{v})| = 1$ then
    return $L(\ta{u}) \cap L(\ta{v})$.}

\item \label{m1b} \texttt{if necessary, interchange left and right
    subtrees in $\ta{u}$ and/or $\tb{v}$ so that $t_{\ell(u)r(v)} +
    t_{r(u)\ell(v)} \leq t_{\ell(u)\ell(v)} + t_{r(u)r(v)}$ and
    $t_{\ell(u)\ell(v)} \leq t_{r(u)r(v)}.$}

\item \label{m1c} \texttt{if {$(t_{\ell(u)\ell(v)} > 0) $} then}
    \begin{enumerate}
    \item \label{m1c1} \texttt{select any leaf $z$ from $L(\ta{\ell(u)})
        \cap L(\ta{\ell(v)})$,}
    \item \label{m1c2} \texttt{return $\{z\} \cup
        \mathtt{Match1}(\ta{r(u)},\tb{r(v)})$.}
    \end{enumerate}
  \item \label{m1d} \texttt{if {$t_{r(u)\ell(v)} = 0$} then return
      $\mathtt{Match1}(\ta{u},\tb{r(v)})$.}

  \item \label{m1e} \texttt{if {$t_{\ell(u)r(v)} = 0$} then return
      $\mathtt{Match1}(\ta{r(u)},\tb{v})$.}

  \item \label{m1f} \texttt{if {$(t_{\ell(u)r(v)}+t_{r(u)\ell(v)} \geq
        \delta t_{uv})$} then}
    \begin{enumerate}
    \item \label{m1f1} \texttt{if necessary, interchange the left and
        the right subtrees of both $\ta{u}$ and $\tb{v}$ so that
        $t_{\ell(u)r(v)} \leq t_{r(u)\ell(v)}$,}
    \item \label{m1f2} \texttt{select any leaf $z$ from $L(\ta{\ell(u)})
        \cap L(\ta{r(v)})$,}
    \item \label{m1f3} \texttt{return $\{z\} \cup
        \mathtt{Match1}(\ta{r(u)},\tb{\ell(v)})$.}
    \end{enumerate}
  \item \label{m1g}\texttt{return
      $\mathtt{Match1}(\ta{r(u)},\tb{r(v)})$.}
\end{enumerate} \medskip

\citet{sw1993} devised an efficient polynomial time algorithm for
finding the maximum agreement subtree of two given binary trees. Our
algorithm is not optimal, but it is easy to analyze; moreover, it is
easy to guarantee a lower bound on the returned value.  We show how the
value of~$\alpha$ in Lemma~\ref{lemma:match1} may also be obtained by
analyzing the algorithm.

\begin{proof}[Alternative proof of Corollary~\ref{cor:match1}.] Set $x
  \coloneqq \rho(T_1)$ and $y \coloneqq \rho(T_2)$. We now analyze the
  execution of the a call to \texttt{Match1}$(T_1^x, T_2^y)$.

  First observe that, each time an instance of algorithm \texttt{Match1}
  is being executed, it calls itself recursively only once in that
  instance. That happens, say $k$ times, each time going deeper in the
  recursion levels, until a base case is reached in
  line~\ref{m1a}. Thus, we may define two sequences of nodes $x = x_0,
  x_1, x_2, \dots, x_k$ and $y = y_0, y_1, y_2, \dots, y_k$ that
  correspond to the roots of the trees passed as arguments in each
  triggered call. That is, \texttt{Match1}$(T_1^{x_0}, T_2^{y_0})$ calls
  \texttt{Match1}$(T_1^{x_1}, T_2^{y_1})$, which in turn calls
  \texttt{Match1}$(T_1^{x_2}, T_2^{y_2})$, and so on.  Note that the
  nodes in each sequence need not be distinct. For example, if
  \texttt{Match1} is called from line~\ref{m1d}, then $x_{i+1} = x_i$;
  and similarly, if \texttt{Match1} is called from line~\ref{m1e}, then
  $y_{i+1} = y_i$.

  As a shorthand notation, define
  $$
  t_i \coloneqq t_{x_i y_i}.
  $$
  
  Now suppose \texttt{Match1}$(T_1^{u}, T_2^{v})$ is being called with
  $u = x_i$ and $v = y_i$. Using our notation, we have $t_i =
  t_{uv}$. For each possibility of calling \texttt{Match1} recursively,
  we obtain a lower bound for $t_{i+1}$ in terms of $t_i$.

  After executing line~\ref{m1b}, we have $t_{r(u)r(v)} \geq
  t_{uv}/4$. Hence, if the recursive call in line~\ref{m1c2} is
  triggered, we have
  \begin{equation}
    \label{eq:3b} t_{i+1} \geq t_i/4.
  \end{equation} If the recursive call in line~\ref{m1d} is triggered,
  it is because $t_{\ell(u)\ell(v)} = 0$ and $t_{r(u)\ell(v)} = 0$, which
  implies $t_{ur(v)} = t_{uv}$. Hence, in that case, we must have
  \begin{equation}
    \label{eq:4n5} t_{i+1} = t_i.
  \end{equation}
  Similarly, if the recursive call in line~\ref{m1e} is triggered, we
  have $t_{r(u)v} = t_{uv}$, which also implies~\eqref{eq:4n5}. The
  conditions in lines~\ref{m1f} and~\ref{m1f1} imply that
  $t_{r(u)\ell(v)} \geq \delta t_{uv}/2$. Hence, if the recursive call
  in line~\ref{m1f3} is made, we must have
  \begin{equation}
    \label{eq:6c} t_{i+1} \geq \delta t_{i}/2.
  \end{equation}
  In line~\ref{m1g}, since $t_{\ell(u)r(v)}+t_{r(u)\ell(v)} < \delta
  t_{uv}$ and $t_{\ell(u)\ell(v)}=0$, we have $t_{r(u)r(v)} >
  (1-\delta)t_{uv}$. Hence, if \texttt{Match1} is called from
  line~\ref{m1g}, then we have
  \begin{equation}
    \label{eq:7} t_{i+1} > (1-\delta)t_{i}.
  \end{equation}

  Now suppose that, during the entire execution of the recursive
  algorithm, line~\ref{m1c2} is executed $a$ times, line~\ref{m1d} is
  executed $b$ times, line~\ref{m1e} is executed $c$ times,
  line~\ref{m1f3} is executed $d$ times, and line~\ref{m1g} is executed
  $e$ times. Since a new leaf $z$ is returned each time one of the
  lines~\ref{m1c2} or~\ref{m1f3} is executed, the set returned by the
  outermost call \texttt{Match1}$(T_1^{x}, T_2^{y})$ (when the execution
  halts) contains precisely $a + d + 1$ leaves. Therefore, it is enough
  to show the following.

  \begin{claim} $a + d > \alpha m$.
  \end{claim}

  To prove the claim consider the sequence $t_0, t_1, \dots, t_k$. For
  each $i$ in $\{0,1,\dots,k-1\}$, we know that $t_{i+1}$ satisfies one
  of~\eqref{eq:3b},\eqref{eq:4n5},\eqref{eq:6c} or
  \eqref{eq:7}. Therefore, we have
  \begin{equation}
    \label{eq:beforefinal} t_{k} \geq t_0 \left(\frac{1}{4}\right)^{a}
    \left(\frac{\delta}{2}\right)^d (1-\delta)^e.
  \end{equation}
  Now observe that $t_0 = 2^m$. Also, since the $k$-th recursive call is
  a base case, we have $t_k = 1$. Moreover, because $\delta < 1/2$, we
  have $1/4^a > (\delta/2)^a$. Using these observations and the fact
  that $e \leq m$, equation~\eqref{eq:beforefinal} yields
  $$
  1 > 2^m \left(\frac{\delta}{2}\right)^{a+d} (1-\delta)^m.
  $$
  Solving for $a + d$, we obtain
  $$
  a + d > \left(\frac{1 + \log(1-\delta)}{1 - \log\delta}\right)\!m =
  \alpha m.
  $$
  The claim is proved, and the corollary follows.
\end{proof}

\subsection {The unrooted case}
\label{sec:unrooted_case}

We now consider the case when one of the trees is unrooted and
balanced. We define two classes of unrooted, balanced trees. The center
of a tree may be either a single vertex or a pair of adjacent
vertices. Let $m \geq 0$ be an integer.  When all leaves of a tree are
at distance $m$ from the center and the center is a single vertex, we
say that the tree is in class~$\mcc_m$. When all leaves are at distance
$m$ from the center and the center is a pair of adjacent vertices, we
say the phylogenetic tree is in class $\mcb_{m + 1}$.  By construction,
trees in the class $\mcc_m$ have $3 \times 2^{m-1}$ leaves, and trees in
the class $\mcb_m$ have $2^{m}$ leaves.

\begin{thm} \label{t1} If $T_1$ is a balanced phylogenetic tree on $n$
  leaves, and $T_2$ is an arbitrary phylogenetic tree on the same
  leaf-set, then they have an agreement subtree on at least $\alpha
  \log{\frac{2n}{3}}$ leaves, where $\alpha$ is the constant defined
  in~\eqref{eq:alpha}.
\end{thm}

\begin{proof} If $T_1$ is in class $\mcb_m$, for some $m$, then $n =
  2^m$. Let $\{x,y\}$ be its central edge. We add a new vertex $z$, and
  replace the edge $\{x,y\}$ by edges $\{x,z\}$ and $\{y,z\}$, and root
  the tree at $z$. For $T_2$, we add a new vertex $w$, replace an
  arbitrary edge $\{u,v\}$ by edges $\{u,w\}$ and $\{v,w\}$, and root
  $T_2$ at $w$. Notice that, for any rooted agreement subtree of $T_1^z$
  and $T_2^w$, we may ignore the root and obtain an unrooted agreement
  subtree of~$T_1$ and~$T_2$. Applying Lemma~\ref{lemma:match1} to
  $\ta{z}$ and $\tb{w}$ gives a lower bound of $\alpha m$ on the size of
  the maximum agreement subtree of $\ta{z}$ and $\tb{w}$. The desired
  bound follows.

  If $T_1$ is in class $\mcc_m$, for some $m$, then $n = 3\times
  2^{m-1}$. Let $z$ be the center of $T_1$. Let $X$ be the set of leaves
  in two of the three branches rooted at $z$. Note that $T_1|X$ is in
  class $\mcb_m$ and $T_2|X$ is an arbitrary phylogenetic
  tree. Proceeding as in the above paragraph, we obtain a lower bound of
  $\alpha m$ on the size of the maximum agreement subtree of $T_1|X$ and
  $T_2|X$. Hence, $T_1$ and $T_2$ have an agreement subtree on at least
  $\alpha \log{\frac{2n}{3}}$ leaves.
\end{proof}

The following proposition for the case when one of the trees is ``almost
balanced'' is proved with little extra effort.

\begin{prop} \label{p2} For every $k > 0$, there is a constant $\alpha_k
  > 0$ such that, if $T_1$ and $T_2$ are binary trees on the same
  leaf-set of cardinality $n$, and $T_1$ has radius at most $k\log n-1$,
  then they have a maximum agreement subtree on at least $\alpha_k \log
  n$ leaves.
\end{prop}

\begin{proof} In tree $T_1$, we subdivide the central edge (if it has a
  central edge) or an edge adjacent to the center (if its center is a
  single vertex), and root the tree at the newly inserted vertex of
  degree 2. (We have bounded the radius of $T_1$ by $k\log n-1$ and not
  $k\log n$ only to allow the possibility that when we root $T_1$, its
  radius may increase by 1.)  We root $T_2$ by subdividing an
  arbitrarily chosen edge. We then construct a rooted, balanced, binary
  tree $T_1^\prime$ of height $k\log n$ that contains $T_1$ as a subtree
  (in the sense that $T_1 \preceq T_1^\prime$). Now by
  Lemma~\ref{lemma:match1}, we assert that $T_1^\prime$ and $T_2$ (hence
  also $T_1$ and $T_2$) have an agreement subtree on at least
  $$
  \frac{k \log{n} \log(1-\delta) + \log{n}}{1-\log{\delta}} = \left(
    \frac{1 + k \log(1-\delta)}{1 - \log \delta} \right) \log n
  $$ 
  leaves. We select $\delta $ sufficiently small to satisfy $1 + k
  \log(1-\delta) > 0$, and set
  \[ 
  \alpha_k \coloneqq \frac{1 + k \log(1-\delta)}{1 - \log \delta}.
  \]
  Then there is an agreement subtree on at least $\alpha_k \log n$
  leaves.  Indeed, the above value of $\alpha_k$ may also be obtained by
  (re)analyzing algorithm \texttt{Match1} as in the alternative proof of
  Corollary~\ref{cor:match1}.
\end{proof}

\section{General binary trees}

Our approach to general binary trees is based on the following
intuition: every binary tree has large diameter or contains (as a
restriction) a balanced subtree of large height.

For $0 \leq k \leq h$, let $f(h,k)$ be the maximum number of leaves a
rooted tree of height at most $h$ can have so that no restriction of the
tree is a balanced, binary tree of height more than $k$.

\begin{lem}
  \label{lemma:fhk} If $h=k$ or $k=0$, then $f(h,k) = 2^k$. If $0 < k <
  h$, then
  \[ f(h,k) = \sum_{i=0}^{k} \binom{h-i-1}{k-i}2^{i}.
  \]
\end{lem}

\begin{proof} We claim the following recurrence for $f(h,k)$:
  \begin{equation}\label{e6} f(h,k) =
    \begin{cases} 2^k & \text{ if } h = k \text{ or } k = 0,\\
      f(h-1,k)+f(h-1,k-1) & \text{ if } 0 < k < h.
    \end{cases}
  \end{equation} 
  This is proved as follows. If $h = k$ or $k = 0$, then we have $f(h,k)
  = 2^k$, the extremal tree being the rooted, balanced tree of height
  $k$.

  Now suppose that $h > k > 0$. We first prove that $f(h,k) \leq
  f(h-1,k)+f(h-1,k-1)$. Let $T$ be a binary tree of height at most $h$
  with more than $f(h-1,k)+f(h-1,k-1)$ leaves. Suppose $T$ has $x$
  leaves in the left subtree and $y$ leaves in the right
  subtree. Without loss of generality assume $y \leq x$. If $x >
  f(h-1,k)$, then the left subtree of $T$ would have a restriction to a
  balanced, binary tree of height $k+1$. Therefore, we may assume that
  $x \leq f(h-1,k)$, which implies $f(h-1,k-1) < y \leq x$. It follows
  that both the left and the right subtrees have restrictions to
  balanced trees of height~$k$, and that $T$ has a restriction to a
  balanced tree of height $k+1$, which is a contradiction.

  Next we show that $f(h,k) \geq f(h-1,k)+f(h-1,k-1)$. Consider the tree
  $T(h,k)$ defined as follows: if $h=k$ or $k=0$, then $T(h,k)$ is a
  balanced, binary tree of height $k$; otherwise, its left subtree is an
  extremal tree for parameters $h-1$ and $k$, and its right subtree is
  an extremal tree for parameters $h-1$ and $k-1$. Thus $T(h,k)$ has
  precisely $f(h-1,k)+f(h-1,k-1)$ leaves, and does not contain a
  restriction that is a balanced tree of height more than $k$.

  Thus we have $f(h,k) = f(h-1,k)+f(h-1,k-1)$ for $0 < k < h$, and the
  tree $T(h,k)$ constructed above is an example of an extremal tree for
  parameters $h$ and $k$. In fact the above arguments, together with
  induction on $h+k$, show that $T(h,k)$ is the unique such tree. We
  skip the details.

  Now the solution to the recurrence relation is obtained by expanding
  it until all terms are expressed as $f(i,i) = 2^i$ for some $i > 0$ or
  $f(j,0) = 1$ for some $j > 0$.
  \begin{figure}[H]
    \begin{center}
      \includegraphics[scale=1]{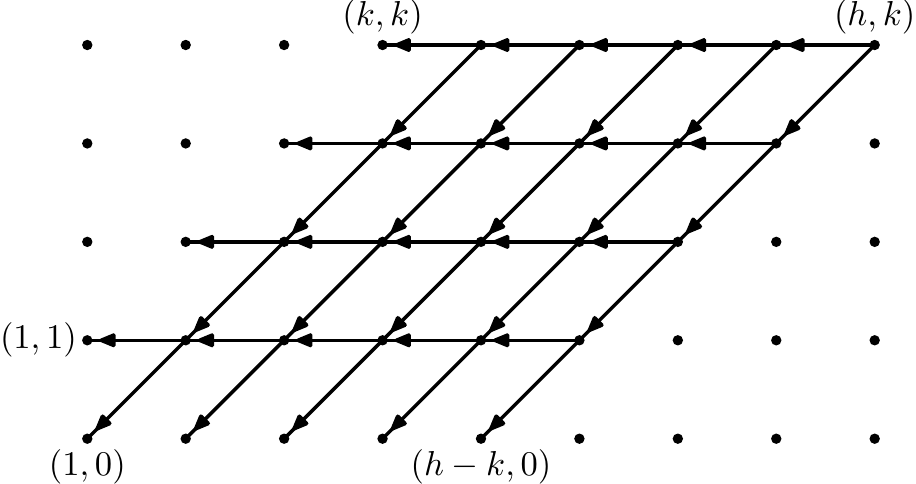}
    \end{center}
    \caption[]{An illustration for the recurrence in~\eqref{e6}}
    \label{fig-trees}
  \end{figure} 
  In Figure~\ref{fig-trees} above, each directed path from the point
  $(h,k)$ to $(i,i)$ contributes the term $f(i,i) = 2^i$, and there are
  $\binom{h-i-1}{k-i}$ such paths; similarly every path from $(h,k)$ to
  $(j,0)$ contributes $f(j,0) = 1$, and there are $\binom{k + j - 1}{j}$
  such paths. Hence, for $0 < k < h$, we have
  $$
  f(h,k) = \sum_{i=1}^{k} \binom{h-i-1}{k-i}2^{i} + \sum_{j=0}^{h-k-1}
  \binom{k+j-1}{j}
  $$
  Since the second sum is $\binom{h-1}{k}$, the desired bound follows.
\end{proof}

\begin{cor}
\label{cor:fhk} For $1 \leq k \leq h$, we have $f(h,k) \leq (2h)^k$.
\end{cor}

\begin{proof} We apply Lemma~\ref{lemma:fhk}. If $k = 1$, then $f(h,k) =
  h + 1 \leq (2h)^k$. If $1 < k < h$, then we have
  $$
  f(h,k) = \sum_{i=0}^{k} \binom{h-i-1}{k-i}2^{i} < {h \choose k}
  \sum_{i=0}^{k} 2^{i} < (2h)^k.
  $$
  If $k = h$, then $f(h,k) = 2^k \leq (2h)^k$.
\end{proof}

Define $\displaystyle \phi(n, a) = \frac{(\log n)^a}{2}$ and
$\displaystyle \psi(n,b) = \frac{(\log n)^b}{\log \log n}$.

\begin{cor}
  \label{cor:Ramsey} Given any $a, b \in (0,1)$ such that $a + b = 1$,
  every tree with~$n > 2$ leaves contains either a path of length at
  least $(\log n)^{\psi(n, b)}$ or a balanced subtree of height at least
  $\phi(n, a)$.
\end{cor}

\begin{proof} Let $a, b \in (0,1)$ such that $a + b = 1$. Let $k \leq
  \phi(n, a)$ and $h \leq (\log n)^{\psi(n, b)}$.  We have
  \begin{eqnarray}
    \label{eq:ramsey} k \log (2h) &=& k + k \log h \nonumber \\ & \leq &
    \phi(n, a) + \phi(n, a)\psi(n, b)\log \log n \nonumber \\ & = &
    \frac{(\log n)^a}{2} + \frac{\log n}{2} \nonumber \\ & < & \log n.
  \end{eqnarray}
  Hence, by Corollary~\ref{cor:fhk}, we conclude that $f(h,k) < n$. Now
  the corollary follows.
\end{proof}

\begin{prop}\label{p4} If $T_1$ and $T_2$ are binary trees on the same
  leaf-set of cardinality $n$, and $T_1$ is a caterpillar, then they
  have a maximum agreement subtree on at least $\frac{1}{3}\log n$
  leaves.
\end{prop}

\begin{proof}
  The proof of this fact goes along the lines of the proof of Theorem~1
  in~\citet{ss2009}. We sketch it here. We embed $T_1$ and $T_2$ in the
  plane so that the leaves of $T_1$ are on one side of the longest path
  in $T_1$. Without loss of generality, suppose that the leaves of $T_1$
  appear in the order $1,2,\ldots, n$. The embedding of $T_2$ imposes a
  circular order on its leaves. We cut this circular order arbitrarily
  to get a linear order $i_1,i_2, \ldots, i_n$. Next we find the longest
  monotone subsequence of $i_1,i_2, \ldots, i_n$; it has length at least
  $\sqrt{n}$ by the Erd\H{o}s-Szekeres Theorem~\cite{es1935}. Let $X$ be
  the set of leaves in this subsequence. We restrict $T_1$ and $T_2$ to
  $X$ obtaining $T_1|X$ and $T_2|X$. Notice that $T_1|X$ is still a
  caterpillar. We further restrict both trees to $Y \subseteq X$ so that
  $T_2|Y$ is a caterpillar with a maximum number of leaves. Thus $|Y|$
  is at least $\log n$ (the extremal case being when $T_2|X$ is
  balanced). Now both $T_1|Y$ and $T_2|Y$ are caterpillars (see
  Figure~\ref{fig:caterpillar_embeddings} below). Now let $y_1, y_2,
  \dots, y_k$ be the elements of $Y$ in the order they appear in the
  embedding of~$T_1$.
  \begin{figure}[H] \centering
    \includegraphics{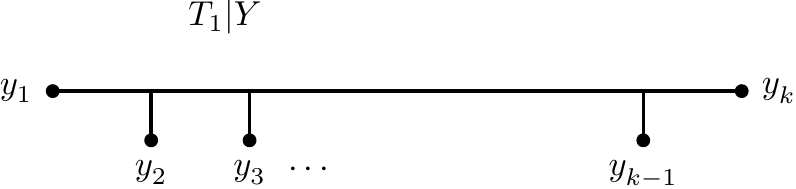} \\ \bigskip
    \includegraphics{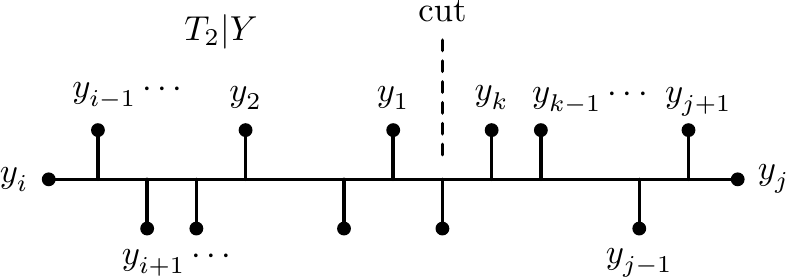}
    \caption{The embeddings of $T_1|Y$ and $T_2|Y$.}
    \label{fig:caterpillar_embeddings}
  \end{figure}
  In Figure~\ref{fig:caterpillar_embeddings}, we can see that there are
  three maximal agreement caterpillars, namely, caterpillars with
  leaf-sets $\{y_1, y_2, \dots, y_i\}$, $\{y_i, y_{i+1}, \dots, y_j\}$
  and $\{y_j, y_{j+1}, \dots, y_k\}$. One of them must have length at
  least $(k + 2)/3 \geq \frac{1}{3} \log n$.
\end{proof}

\begin{thm}
  \label{thm:root_log} If $T_1$ and $T_2$ are binary trees on the same
  leaf-set of cardinality~$n > 2$, then they have a maximum agreement
  subtree having at least $\frac{\alpha}{2} \sqrt{\log n} + \alpha \log
  \frac{2}{3}$ leaves.
\end{thm}

\begin{proof} Applying Corollary~\ref{cor:Ramsey} (with $a = b = 1/2$),
  one of the trees must contain a balanced subtree of height at least
  $\phi(n,1/2)$ or a path of length at least $(\log
  n)^{\psi(n,1/2)}$. Suppose that one of the trees contains a balanced
  subtree of height at least $\phi(n,1/2)$. Let $A$ be the leaf-set of
  such a balanced subtree. Therefore, after restricting the other tree
  to $A$, we can claim by Theorem~\ref{t1} that $T_1$ and $T_2$ have a
  common subtree on at least $\alpha \phi(n,1/2) + \alpha \log
  \frac{2}{3}$ leaves from~$A$.

  If such a balanced subtree does not exist in either of the two trees,
  then there is a path (and hence a caterpillar) of length at least
  $(\log n)^{\psi(n,1/2)}$ in one of the trees. We restrict both trees
  to the set of leaves in this caterpillar. Therefore, by
  Proposition~\ref{p4}, there must be a common subtree on at least
  $\frac{1}{3}\psi(n,1/2)\log \log n = \frac{1}{3}\sqrt{\log n}$
  leaves. Taking the maximum value of~$\alpha$ as in
  Remark~\ref{rem:alpha}, this is a quantity larger than the desired
  bound.
\end{proof}

\section{When both trees are balanced}

We now investigate the size of a maximum agreement subtree of two
balanced, binary trees. In this case, in Theorem~\ref{t2} we obtain a
much better bound than that of Theorem~\ref{t1}.

\begin{lem}
  \label{lemma:match2} Suppose $T_1$ and $T_2$ are rooted, balanced,
  binary trees of height $m_1$ and $m_2$, respectively. Suppose that
  $|L(T_1) \cap L(T_2)| = t > 0$. Then for all $\delta \in
  \left(0,\frac{1}{4}\right)$, the two trees have a rooted maximum
  agreement subtree on at least $2^{g(m_1,m_2,t)}$ leaves, where
  $$
  g(m_1,m_2,t) \coloneqq \frac{(m_1+m_2)\log(1-3\delta) +
    \log{t}}{\log{(1-3\delta)}-\log{\delta}}.
  $$
\end{lem}

\begin{proof} Let $M(m_1,m_2,t)$ be the minimum value of
  $\mast\{T_1,T_2\}$ (over all choices of $T_1$ and $T_2$), where $T_1$
  and $T_2$ are as in the statement of the lemma. Observe that
  $M(m_1,m_2,t)$ is a monotonically non-decreasing function of $t$. We
  show the result by induction on $m_1+m_2$.

  \

  \noindent {\em Base case:} When $m_1 + m_2 \in \{0,1\}$, at least one
  of the trees has a single vertex (which is its leaf and root), hence
  $\mast\{T_1,T_2\} = 1$. Also, since $t = 1$, we have $g(m_1,m_2,t)
  \leq 0$, and the claim is true. So we assume below that $m_1 \geq 1$
  and $m_2 \geq 1$.

  \

  \noindent {\em Induction step:} Let $u \coloneqq \rho(T_1)$ and $v
  \coloneqq \rho(T_2)$. As in Theorem~\ref{t1}, we have $t = t_{uv} =
  t_{\ell(u)\ell(v)} + t_{r(u)r(v)} + t_{\ell(u)r(v)} + t_{r(u)\ell(v)}$
  and we assume, without loss of generality, that $t_{\ell(u)\ell(v)} +
  t_{r(u)r(v)} \geq t_{\ell(u)r(v)} + t_{r(u)\ell(v)}$ and $t_{r(u)r(v)}
  \geq t_{\ell(u)\ell(v)}$ to obtain
  \begin{equation}
    \label{eq:case1a_again} t_{r(u)r(v)} \geq \lceil t/4 \rceil.
  \end{equation}

  \noindent {\em Case 1:} $t_{\ell(u)\ell(v)} \geq \delta t$.

  By~\eqref{eq:case1a_again}, we also have $t_{r(u)r(v)} \geq \delta
  t$. In this case, we take a maximum agreement subtree $S_\ell$ of
  $\ta{\ell(u)}$ and $\tb{\ell(v)}$, and a maximum agreement subtree
  $S_r$ of $\ta{r(u)}$ and $\tb{r(v)}$. We then construct $S = S_\ell
  \circ S_r$, which by Proposition~\ref{lemma:join} is an agreement
  subtree of $T_1$ and $T_2$. Therefore, we have
  \begin{eqnarray*}
    \label{eq:case1_match2_a} M(m_1,m_2,t) & \geq & 2M(m_1-1,m_2-1,
    \lceil \delta t \rceil)\nonumber \\ & \geq & 2\times
    2^{g(m_1-1,m_2-1,\lceil \delta t \rceil)} \nonumber \\ & \geq &
    2^{(1+g(m_1-1,m_2-1,\lceil \delta t \rceil))} \nonumber \\ & \geq &
    2^{g(m_1,m_2, t)},
  \end{eqnarray*} where the last step follows from
  \begin{eqnarray*}
    \label{eq:case1_match2_b} 1+g(m_1-1,m_2-1,\lceil \delta t \rceil)
    &=& 1 + \frac{(m_1-1+m_2-1)\log(1-3\delta) + \log (\delta t)}
    {\log(1-3\delta)-\log\delta } \nonumber \\ &=&
    \frac{(m_1+m_2)\log(1-3\delta) + \log t - \log(1-3\delta)}
    {\log(1-3\delta)-\log\delta } \nonumber \\ & > & g(m_1,m_2,t),
  \end{eqnarray*} where the last inequality requires that $\delta \in
  \left(0,\frac{1}{4}\right)$.

  \

  \noindent {\em Case 2:} $t_{\ell(u)r(v)} \geq \delta t$ and
  $t_{r(u)\ell(v)} \geq \delta t$.

  The calculation in this case is identical to that of Case 1, so we
  omit it.

  \

  \noindent {\em Case 3:} $t_{\ell(u)r(v)} < \delta t$ and
  $t_{r(u)\ell(v)} < \delta t$.

  Since Case 1 has been examined, we assume that $t_{\ell(u)\ell(v)} <
  \delta t$, which implies $t_{r(u)r(v)} > (1 -3\delta) t$. Since
  $\mast\{T_1,T_2\}$ must be at least $\mast\{\ta{r(u)},\tb{r(v)}\}$, we
  have $M(m_1,m_2,t) \geq M(m_1-1,m_2-1, \lceil (1-3\delta)
  t)\rceil$. Now the result follows from the assumption that $\delta \in
  \left(0,\frac{1}{4}\right)$ and the following:
  \begin{eqnarray*}
    \label{eq:case3_match2} && g(m_1-1,m_2-1,\lceil (1-3\delta)t \rceil)
    \nonumber \\ & \geq & \frac{(m_1-1+m_2-1)\log(1-3\delta) +
      \log(1-3\delta)t} {\log(1-3\delta)-\log\delta} \nonumber \\ & = &
    \frac{(m_1+m_2)\log(1-3\delta) + \log t - \log(1-3\delta)}
    {\log(1-3\delta)-\log\delta} \nonumber \\ & > & g(m_1,m_2,t).
  \end{eqnarray*}

  \

  \noindent {\em Case 4:} $t_{\ell(u)r(v)} < \delta t$ and
  $t_{r(u)\ell(v)} \geq \delta t$.
  
  Since Case 1 has been examined, we assume that $t_{\ell(u)\ell(v)} <
  \delta t$, which implies $t_{r(u)\ell(v)} + t_{r(u)r(v)} > (1-2\delta)
  t$. In this case, since $\mast\{T_1,T_2\}$ must be at least
  $\mast\{\ta{r(u)},T_2\}$, we can write $ M(m_1,m_2,t) \geq M(m_1-1,
  m_2, \lceil (1-2\delta) t \rceil).$ Now the result follows from the
  assumption that $\delta \in \left(0,\frac{1}{4} \right)$ and the
  following:
  \begin{eqnarray*}
    \label{eq:case4_match2} && g(m_1-1, m_2, \lceil (1-2\delta) t\rceil)
    \nonumber \\ & \geq & \frac{(m_1-1+m_2)\log(1-3\delta) +
      \log(1-2\delta)t} {\log(1-3\delta)-\log\delta}\nonumber \\ & = &
    \frac{(m_1+m_2)\log(1-3\delta) + \log t + \log(1-2\delta) - \log
      (1-3\delta)} {\log(1-3\delta)-\log\delta}\nonumber \\ & > &
    g(m_1,m_2,t).
  \end{eqnarray*}

  \

  \noindent {\em Case 5:} $t_{\ell(u)r(v)} \geq \delta t$ and
  $t_{r(u)\ell(v)} < \delta t$.

  The analysis of this case is similar to Case 4, except that we have
  the inequality $M(m_1,m_2,t) \geq M(m_1,m_2-1, \lceil (1-2\delta)
  t\rceil)$.
\end{proof}

\begin{cor}
  \label{cor:match2} Let $$\delta \in
  \left(0,\frac{1}{3}-\frac{1}{3\sqrt{2}}\right) \text{~~and~~} \beta
  \coloneqq \left(\frac{1 +2\log (1 - 3\delta)}{\log (1-3\delta) - \log
      \delta}\right).
  $$ 
  If $T_1$ and $T_2$ are rooted, balanced, binary trees on the same
  leaf-set of cardinality $2^m$, then $T_1$ and $T_2$ have a maximum
  agreement subtree on at least $2^{\beta m}$ leaves.
\end{cor}

\begin{proof} We set $m_1 = m_2 = m$ and $t = 2^m$ in
  Lemma~\ref{lemma:match2}. Moreover, we now require $\delta $ to be
  less that $\left(\frac{1}{3}-\frac{1}{3\sqrt{2}}\right)$ (which is
  less that 1/4) so as to ensure that $\beta$ is positive.
\end{proof}

As in Section~\ref{sec:onebalanced}, we present algorithm
\texttt{Match2} that closely follows the recursions in the ṕroof of
Lemma~\ref{lemma:match2}. It takes as input two rooted, balanced, binary
trees, and returns a set of leaves in a common subtree. Algorithm
\texttt{Match2} depends on a real positive $\delta $, which we require
to be sufficiently small for the algorithm to give a desired bound on
the size of a common subtree. The algorithm is somewhat greedy and
suboptimal. The analysis of the performance of \texttt{Match2} makes
Lemma~\ref{lemma:match2} much more transparent, giving an alternative
proof of Corollary~\ref{cor:match2}. We then apply the corollary to
prove the main results of this section for unrooted, balanced (or
``almost balanced'') trees.

\newpage

\noindent \texttt{Algorithm Match2($\ta{u}$, $\tb{v}$)}

\begin{enumerate}[1:]

\item \label{m2a} \texttt{if $|L(\ta{u})| = 1$ or $|L(\ta{v})| = 1$ then
    \\ return $L(\ta{u}) \cap L(\ta{v})$.}

\item \label{m2b} \texttt{if necessary, interchange left and right
    subtrees in $\ta{u}$ and/or $\tb{v}$ so that $t_{\ell(u)r(v)} +
    t_{r(u)\ell(v)} \leq t_{\ell(u)\ell(v)} + t_{r(u)r(v)}$ and
    $t_{\ell(u)\ell(v)} \leq t_{r(u)r(v)}$.}
  
\item \label{m3c} \texttt{if $(t_{\ell(u)\ell(v)} \geq \delta t_{uv}$
    and $t_{r(u)r(v)} \geq \delta t_{uv})$ then \\ return
    $\mathtt{Match2}(\ta{\ell(u)},\tb{\ell(v)}) \cup
    \mathtt{Match2}(\ta{r(u)},\tb{r(v)})$.}
  
\item \label{m3d} \texttt{if $(t_{\ell(u)r(v)} \geq \delta t_{uv}$ and
    $t_{r(u)\ell(v)} \geq \delta t_{uv})$ then \\ return
    $\mathtt{Match2}(\ta{\ell(u)},\tb{r(v)}) \cup
    \mathtt{Match2}(\ta{r(u)},\tb{\ell(v)})$.}
  
\item \label{m3e} \texttt{if $(t_{\ell(u)r(v)} < \delta t_{uv}$ and
    $t_{r(u)\ell(v)} < \delta t_{uv})$ then \\ return
    $\mathtt{Match2}(\ta{r(u)},\tb{r(v)})$.}
  
\item \label{m3f} \texttt{if $(t_{\ell(u)r(v)} < \delta t_{uv}$ and
    $t_{r(u)\ell(v)} \geq \delta t_{uv})$ then \\ return
    $\mathtt{Match2}(\ta{r(u)},\tb{v})$.}

\item \label{m3g} \texttt{if $(t_{\ell(u)r(v)} \geq \delta t_{uv}$ and
    $t_{r(u)\ell(v)} < \delta t_{uv})$ then \\ return
    $\mathtt{Match2}(\ta{u},\tb{r(v)})$.}

\end{enumerate}

We now analyze the above algorithm to compute $\beta $ in
Corollary~\ref{cor:match2} more transparently.

\begin{proof}[Alternative proof of Corollary~\ref{cor:match2}.]
  We prove the result by analyzing \texttt{Match2}. In the beginning, we
  call $\mathtt{Match2}(\ta{x},\tb{y})$, where $x \coloneqq \rho(T_1)$
  and $y \coloneqq \rho(T_2)$.  Let $\mct$ be the tree of recursive
  calls to \texttt{Match2} constructed as follows: the pair $(x,y)$ is
  the root of $\mct$. If $\texttt{Match2}(\ta{u},\tb{v})$ is called
  during the execution of the algorithm, then $(u,v)$ is a vertex of
  $\mct$. If $\texttt{Match2}(\ta{u},\tb{v})$ calls
  $\texttt{Match2}(\ta{u'},\tb{v'})$, then $(u',v')$ is a child of
  $(u,v)$. The leaf vertices of $\mct$ correspond to the function calls
  that return in line~\ref{m2a}. Observe that in line~\ref{m2a}, a set
  containing a single new leaf is returned. By construction, the number
  of leaves in the common subtree returned by \texttt{Match2} is
  precisely the number of leaves of $\mct$.

  The ideas in this lemma are similar to those in the alternative proof
  of Corollary~\ref{cor:match1}. We consider an arbitrary root-to-leaf
  path in $\mct$ and we show that it branches at least $\beta m$
  times. We then conclude that $\mct$ has at least $2^{\beta m}$ leaves,
  thereby proving the theorem.

  Now consider an arbitrary root-to-leaf path $(x_0, y_0), (x_1, y_1),
  \dots, (x_k, y_k)$ in~$\mct$, where $(x_0,y_0) \coloneqq (x,y)$.  As a
  shorthand notation, define
  $$
  t_i \coloneqq t_{x_i y_i}.
  $$
  
  Now suppose \texttt{Match2}$(T_1^{u}, T_2^{v})$ is being called with
  $u = x_i$ and $v = y_i$. Using our notation, we have $t_i =
  t_{uv}$. For each possibility of calling \texttt{Match2} recursively,
  we obtain a lower bound for $t_{i+1}$ in terms of $t_i$.

  First observe that, as in the case of \texttt{Match1}, we relabel
  $\ell(u),r(u)$ and $\ell(v),r(v)$ so that we have $t_{\ell(u)\ell(v)}
  + t_{r(u)r(v)} \geq t_{uv}/2$ and $t_{r(u)r(v)} \geq
  t_{\ell(u)\ell(v)}$ (which implies $t_{r(u)r(v)} \geq
  t_{uv}/4$). Hence, for the choice of $\delta$, we have
  \begin{equation} t_{r(u)r(v)} \geq t_{uv}/4 \geq \delta
    t_{uv}.\label{eq:one_fourth}
  \end{equation}
  If a recursive call in line~\ref{m3c} is triggered, then $(x_{i+1},
  y_{i+1})$ is either $(\ell(u),\ell(v))$ or $(r(u),r(v))$. In both
  cases, under the conditions in line~\ref{m3c}, we have
  \begin{equation}
    \label{eq:line3recur} t_{i+1} \geq \delta t_i.
  \end{equation}
  Similarly, if a recursive call in line~\ref{m3d} is triggered, then we
  also have~\eqref{eq:line3recur}.

  After lines~\ref{m3c} and~\ref{m3d}, we must have $t_{\ell(u) \ell(v)}
  < \delta t_{uv}$ because $t_{\ell(u) \ell(v)} \leq t_{r(u) r(v)}$ and
  the condition in line~\ref{m3c} has failed.
  
  If the recursive call in line~\ref{m3e} is triggered, then $(x_{i+1},
  y_{i+1}) = (r(u), r(v))$, and
  \begin{equation}
    \label{eq:line5recur} t_{i+1} \geq (1 - 3\delta) t_i,
  \end{equation}
  because $t_{\ell(u)\ell(v)} < \delta t_{uv}$, $t_{\ell(u)r(v)} <
  \delta t_{uv}$ and $t_{r(u)\ell(v)} < \delta t_{uv}$.

  If the recursive call in line~\ref{m3f} is triggered, then $(x_{i+1},
  y_{i+1}) = (r(u), v)$, and
  \begin{equation}
    \label{eq:line6recur} t_{i+1} \geq (1 - 2\delta) t_i,
  \end{equation}
  because $t_{\ell(u)\ell(v)} < \delta t_{uv}$ and $t_{\ell(u)r(v)} <
  \delta t_{uv}$.

  Similarly, if the recursive call in line~\ref{m3g} is triggered, then
  $(x_{i+1}, y_{i+1}) = (u, r(v))$, and~\eqref{eq:line6recur} holds
  because $t_{\ell(u)\ell(v)} < \delta t_{uv}$ and $t_{r(u)\ell(v)} <
  \delta t_{uv}$.

  Along the chosen path in $\mct$, suppose that line~\ref{m3c} is
  executed $a$ times, line~\ref{m3d} is executed $b$ times,
  line~\ref{m3e} is executed $c$ times, line~\ref{m3f} is executed $d$
  times, and line~\ref{m3g} is executed $e$ times. In each of the
  recursive calls, the height of one of the trees decreases by
  1. Therefore, $a + b + c + d + e \leq 2m$. Consequently, we have
  \begin{eqnarray}
    \label{eq:match2} t_{uv} & \geq & t_0
    \delta^{a+b}(1-3\delta)^c(1-2\delta)^{d+e} \nonumber \\ & > & 2^m
    \delta^{a+b}(1-3\delta)^{c + d + e} \nonumber \\ & \geq & 2^m
    \delta^{a+b}(1-3\delta)^{2m-a-b} \nonumber \\ & = & 2^m
    \left(\frac{\delta}{1-3\delta}\right)^{a+b}(1-3\delta)^{2m}
  \end{eqnarray} 
  Since $(x_k, y_k)$ is a leaf of $\mct$, we must have $t_{k} =
  1$. Hence the right-hand-side of~(\ref{eq:match2}) must be less than
  $1$, which implies, for the choice of $\beta$, that $a+b > \beta m$.
  Now, a positive $\delta$ less than $\frac{1}{3} - \frac{1}{3\sqrt{2}}$
  guarantees that $\beta$ is positive.

  We have shown that each root-to-leaf path in~$\mct$ branches at least
  $a+b \geq \beta m$ times, which further implies that there must be at
  least $2^{\beta m}$ leaves in $\mct$. Hence $T_1$ and~$T_2$ must have
  a common subtree on at least $2^{\beta m}$ leaves.
\end{proof}

\begin{rem} Observe that in the above analysis, we showed that each
  root-to-leaf path in $\mct$ has length at least $\beta m$, which
  implies that it is possible to find an agreement subtree (not
  necessarily a {\em maximum} agreement subtree) with at least
  $2^{\beta m}$ leaves that is also balanced (and of height at least
  $\beta m$). The agreement subtree obtained by algorithm
  \texttt{Match2} is illustrated in
  Figure~\ref{fig:balanced-agreement}. We may choose a single leaf
  from each subtree rooted at depth $\lceil \beta m \rceil$, and
  restrict the tree to chosen leaves to obtain a balanced agreement
  subtree with precisely $2^{\lceil \beta m \rceil}$ leaves.
  \begin{figure}[H]
    \centering
    \includegraphics{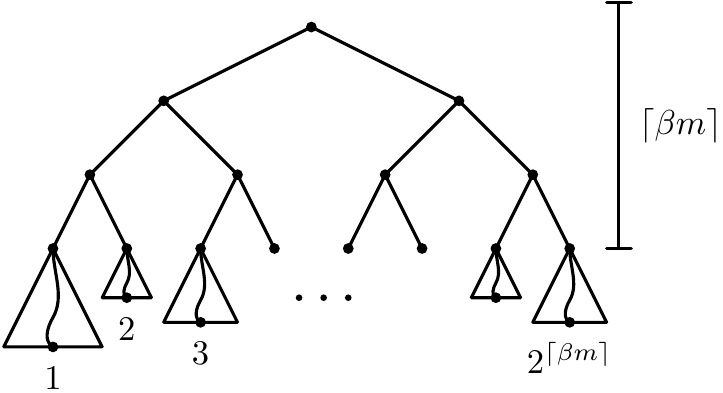}
    \caption{}
    \label{fig:balanced-agreement}
  \end{figure}
  One of the implications of this observation, which we state without
  proof, is that Lemma~\ref{lemma:match2} together with algorithm
  \texttt{Match2} may be used to obtain a lower bound of $2^{\gamma
    m}$ for an agreement subtree of more than 2 balanced binary trees
  of height $m$ for a sufficiently small positive $\gamma$. For
  example, we call algorithm \texttt{Match2} for two rooted, balanced
  trees $T_1$ and $T_2$. The resulting agreement subtree contains a
  rooted, balanced agreement subtree, say $T_{12}$, of height at least
  $\beta m$. We then call algorithm \texttt{Match2} for $T_{12}$ and
  $T_3$. The resulting agreement subtree contains a rooted, balanced
  agreement subtree $T_{123}$, and so on.
\end{rem}

\begin{thm} \label{t2} There exists a constant $c > 0$ such that, if
  $T_1$ and $T_2$ are balanced, binary trees on the same leaf-set, both
  in $\mcb_m$ or both in $\mcc_m$, then they have a maximum agreement
  subtree on at least $2^{\beta m - c}$ leaves.
\end{thm}

\begin{proof}
  As in Theorem~\ref{t1}, we consider the two cases: the trees are
  either both in class $\mcb_m$ or both in class $\mcc_m$. When the
  trees are both in class $\mcb_m$, the proof is analogous to the
  corresponding case in Theorem~\ref{t1}, except that it invokes
  Corollary~\ref{cor:match2} instead of Lemma~\ref{lemma:match1}.

  When the trees are in class $\mcc_m$, the analysis differs only
  slightly from that in Theorem~\ref{t1}. We delete one of the branches
  of $T_1$ rooted at the center, and root the resulting tree at the
  center (which now has degree 2). Let $X$ be the leaf-set of the pruned
  tree. We cannot simply take a restriction of $T_2$ to the leaf-set $X$
  as in Theorem~\ref{t1}, since $T_2|X$ may not be a balanced tree. We
  instead delete one of the branches of $T_2$ rooted at its center, and
  root the pruned tree at its center. Let $Y$ be the leaf-set of the
  pruned tree. Now we can ensure that $|X\cap Y| \geq 2^{m+1}/3$ by
  appropriately choosing the branches of $T_1$ and $T_2$ to be
  deleted. We apply \texttt{Match2} to the rooted trees $T_1|X$ and
  $T_2|Y$. The analysis of \texttt{Match2} does not change, except that
  we now have a constant factor 2/3 on the right-hand-side
  of~(\ref{eq:match2}). Therefore, with $a,b,c,d,e$ defined as in the
  alternative proof of Corollary~\ref{cor:match2}, we have
  $$ 
  a+b > m\left(\frac{1 +2\log (1 - 3\delta) -(1/m) \log (3/2)}{\log
      (1-3\delta) - \log \delta}\right).
  $$
  Taking $c = (\log 3 - 1)/(\log (1-3\delta) - \log \delta)$, there are
  at least $2^{\beta m - c}$ leaves in a common subtree.
\end{proof}

In fact, we have a similar result when the two trees are ``almost
balanced''.

\begin{prop} \label{p3} For every $k > 0$, there is a constant $\beta_k
  > 0$ such that, if $T_1$ and $T_2$ are binary trees on a leaf-set of
  cardinality $n$, each of radius at most $k \log n$, then they have a
  maximum agreement subtree on at least $n^{\beta_k}$ leaves.
\end{prop}

\begin{proof}
  The proof is analogous to the alternative proof of
  Corollary~\ref{cor:match2}: the only change is that now we have $a + b
  + c + d + e \leq 2k\log n$. We use a value of $\delta > 0$ such
  that $$\displaystyle \beta_k \coloneqq \frac{1 + 2k\log (1 -
    3\delta)}{\log(1-3\delta)- \log \delta}$$ is positive, and we have
  $a+b \geq \beta_k \log n$.
\end{proof}

Concerning the maximum agreement subtree problem for balanced trees, we
believe in the following.

\begin{conj}
  Any two balanced, rooted, binary trees of height $m$ have an agreement
  subtree on at least~$2^{m/2}$ leaves.
\end{conj}

We now describe an example of a pair of rooted, balanced, binary trees
of height $2k$, for each $k > 0$, which we believe is an extremal
example. Let $T_1$ and $T_2$ be balanced, binary trees of height $2k$,
rooted at $u$ and $v$, respectively, and both drawn top-down. Let the
leaves of~$T_1$ be labelled $1,2,\cdots, 2^{2k}$ from left to right. We
label the leaves of~$T_2$ from left to right according to the sequence
$\mathtt{swap}(1,2,\cdots,2^{2k})$, which we define recursively as
follows:

\begin{enumerate}
\item If $S$ is a sequence of length 1, then
  \[ \mathtt{swap}(S) = S.
  \]
\item If $S$ is a sequence of length $4^i$, with $i > 0$, written as $S
  \coloneqq S_1:S_2:S_3:S_4$ as a concatenation of 4 sequences of length
  $4^{i-1}$ each, then
  \[ \mathtt{swap}(S) \coloneqq
  \mathtt{swap}(S_1):\mathtt{swap}(S_3):\mathtt{swap}(S_2):\mathtt{swap}(S_4).
  \]
\end{enumerate}

\begin{prop}
  \label{prop:rooted_ex} Trees $T_1$ and~$T_2$ have no rooted agreement
  subtree with more than $2^{k}$ leaves.
\end{prop}

\begin{proof}
  We prove the result by induction on $k$. When $k=1$, the trees have 4
  leaves, with the leaves of $T_1$ labelled 1, 2, 3, 4 from left to
  right, while the leaves of $T_2$ labelled 1, 3, 2, 4 from left to
  right. In this case, a rooted agreement subtree cannot have more than
  2 leaves.

  In the general case, the inductive argument goes as follows. Let $u_1$
  and $u_2$ be the children of $\ell(u)$, and let $u_3$ and $u_4$ be the
  children of $r(u)$. Similarly, in $T_2$, we label the grandchildren of
  $v$ by $v_1,v_2,v_3,v_4$. By construction, we have $L(\ta{u_1}) =
  L(\tb{v_1})$, $L(\ta{u_2}) = L(\tb{v_3})$, $L(\ta{u_3}) =
  L(\tb{v_2})$, and $L(\ta{u_4}) = L(\tb{v_4})$. But a rooted agreement
  subtree cannot have leaves from more than two of the sets
  $L(\ta{u_i}), i \in \{1,2,3,4\}$. Therefore, $\mast\{T_1,T_2\} \leq 2
  \times \mast\{\ta{u_1},\tb{v_1}\} \leq 2 \times 2^{k-1} = 2^k$.
\end{proof}

An analogous but slightly weaker statement holds for unrooted trees. We
use the same labelling scheme as in the rooted case, but remove the
roots, i.e., we delete the vertex $u$, and the edges $\{u,\ell(u)\}$ and
$\{u,r(u)\}$, and add an edge $\{\ell(u),r(u)\}$, and similarly make
$T_2$ unrooted.

\begin{prop}
  \label{prop:unrooted_ex} Trees $T_1$ and~$T_2$ have no agreement
  subtree with more than $3\times 2^{k-1}$ leaves.
\end{prop}

\begin{proof}
  We prove the result by induction on $k$. When $k=1$, the trees have an
  agreement subtree on 3 leaves, but not 4. In the general case, as in
  the rooted case, $L(\ta{u_1}) = L(\tb{v_1})$, $L(\ta{u_2}) =
  L(\tb{v_3})$, $L(\ta{u_3}) = L(\tb{v_2})$, and $L(\ta{u_4}) =
  L(\tb{v_4})$. For $i \in \{1,2,3,4\}$, let $A_i $ denote the set of
  leaves from $L(\ta{u_i})$ that are in a maximum agreement subtree
  $R$. But an agreement subtree cannot have leaves from all four sets
  $L(\ta{u_i}), i \in \{1,2,3,4\}$. Therefore, $A_i = \emptyset $ for
  some $i \in \{1,2,3,4\}$.

  \

  \noindent {\em Case 1:} Two or three of the four sets $A_i$ are
  non-empty. Without loss of generality, let $A_1$ and $A_2$ (and
  possibly also $A_3$) be non-empty. For $i \in \{1,2,3\}$, let $x_i$ be
  the most recent common ancestor of leaves in $A_i$ in the rooted
  subtree $\ta{u_i}$. Now we observe that $R|A_1$ rooted at $x_1$ is a
  rooted agreement subtree of $\ta{u_1}$ and $\ta{v_1}$, and $R|A_2$
  rooted at $x_2$ is a rooted agreement subtree of $\ta{u_2}$ and
  $\ta{v_3}$, and, if $A_3$ is non-empty, $R|A_3$ rooted at $x_3$ is a
  rooted agreement subtree of $\ta{u_3}$ and $\ta{v_2}$. Therefore,
  $|L(R)| = |A_1\cup A_2 \cup A_3| \leq 3 \times 2^{k-1}$ (by
  Proposition~\ref{prop:rooted_ex}).

  \

  \noindent {\em Case 2:} A maximum agreement subtree $R$ has leaves
  from only one of the four sets $L(\ta{u_i}), i \in
  \{1,2,3,4\}$. Without loss of generality, let $L(R) \subseteq
  L(\ta{u_1})$. In this case, by induction, $\mast\{T_1,T_2\} \leq 3
  \times 2^{k-2} < 3 \times 2^{k-1}$.
\end{proof}

\subsection*{Acknowledgements} The first author is partially supported
by the grant CNPq Processo 475064/2010-0 and the second author is
supported by the grant CNPq Processo 151782/2010-5. We would like to
thank CNPq, Brazil for supporting our research. We would also like to
thank the support from the project MaCLinC at Universidade de S\~{a}o
Paulo (USP). Finally, we thank the referee for many useful suggestions
that helped us improve our presentation as well as the main result.

\end{document}